\documentclass[english]{article}
\usepackage[T1]{fontenc}
\usepackage[latin9]{inputenc}
\usepackage[a4paper]{geometry}
\usepackage{babel}
\usepackage{calc}
\usepackage{amsmath}
\usepackage{amsthm}
\usepackage{amssymb}
\usepackage{graphicx}
\usepackage[unicode=true,pdfusetitle,
 bookmarks=true,bookmarksnumbered=false,bookmarksopen=false,
 breaklinks=false,pdfborder={0 0 1},backref=false,colorlinks=false]
 {hyperref}

\usepackage{xcolor}
\usepackage[normalem]{ulem} 

\makeatletter
\theoremstyle{plain}
\newtheorem{thm}{\protect\theoremname}
\theoremstyle{plain}
\newtheorem{defn}[thm]{\protect\definitionname}
\theoremstyle{plain}
\newtheorem*{rem*}{\protect\remarkname}
\theoremstyle{plain}
\newtheorem{lem}[thm]{\protect\lemmaname}
\theoremstyle{remark}
\newtheorem*{acknowledgement*}{\protect\acknowledgementname}

\@ifundefined{date}{}{\date{}}
\usepackage{dsfont}
\usepackage{bbm}

\newcommand\blfootnote[1]{%
  \begingroup
  \renewcommand\thefootnote{}\footnote{#1}%
  \addtocounter{footnote}{-1}%
  \endgroup
}

\makeatother

\usepackage[style=alphabetic,maxbibnames=9]{biblatex}
\providecommand{\acknowledgementname}{Acknowledgement}
\providecommand{\definitionname}{Definition}
\providecommand{\lemmaname}{Lemma}
\providecommand{\remarkname}{Remark}
\providecommand{\theoremname}{Theorem}

\begin{document}
\global\long\def\P{\mathbb{P}}%
\global\long\def\E{\mathbb{E}}%
\global\long\def\R{\mathbb{R}}%
\global\long\def\N{\mathbb{N}}%
\global\long\def\K{\mathbb{K}}%
\global\long\def\Z{\mathbb{Z}}%
\global\long\def\del{\partial}%
\global\long\def\dx{\,\mathrm{d}}%

\noindent 
\global\long\def\I{\mathds{1}}%
\global\long\def\eps{\varepsilon}%
\global\long\def\ue{u^\eps}%
\global\long\def\e{\mathrm{e}}%
\global\long\def\C{\mathcal{C}}%
\global\long\def\cG{\mathcal{G}}%
\global\long\def\A{\mathcal{A}}%
\global\long\def\S{\mathcal{S}}%
\global\long\def\cR{\mathcal{R}}%
\global\long\def\fd{\mathfrak{d}}%

\global\long\def\L{\mathbb{L}}%

\title{An ergodic and isotropic zero-conductance model with arbitrarily strong local connectivity}
\author{Martin Heida\thanks{Weierstrass Institute Berlin, Mohrenstr. 39, 10117 Berlin, Germany,  \texttt{martin.heida@wias-berlin.de}}
\qquad{}Benedikt Jahnel \thanks{Technische Universit\"at Braunschweig, Universit\"atsplatz 2, 38106 Braunschweig, Germany, and \\
Weierstrass Institute Berlin, Mohrenstr. 39, 10117 Berlin, Germany,  \texttt{benedikt.jahnel@tu-braunschweig.de}} 
\qquad{}Anh Duc Vu \thanks{Weierstrass Institute Berlin, Mohrenstr. 39, 10117 Berlin, Germany,  \texttt{anhduc.vu@wias-berlin.de}}}
\date{\today}
\maketitle
\begin{abstract}
We exhibit a percolating ergodic and isotropic lattice model in all but at least two dimensions that has zero effective conductivity in all spatial directions and for all non-trivial choices of the connectivity parameter. The model is based on
the so-called randomly stretched lattice where we additionally elongate layers containing few open edges. 
\end{abstract}

\blfootnote{{\bf Keywords:} homogenisation, percolation,  randomly stretched lattice, conductivity} 
\blfootnote{{\bf MSC2020:} primary: 60K35, 60K37; secondary: 60G55 90B18}

\section{Introduction}
Consider a stationary and ergodic model for a randomly perforated material $G\subset\R^{d}$. Examples we have in mind are the supercritical cluster of the Boolean model based on a Poisson point process or its complement. 
From a physical point of view, $G$ could be the perforations of a sponge-like material which allows for the diffusion of some chemicals while the complement would block these chemicals.

We assume that the scale of the perforations is small compared to the macroscopic dimensions of the material, which we express as $G^\eps := \eps G$ and we are interested of the effective conductivity of $G^\eps$ as $\eps\to0$. 
To be more precises, given a bounded domain $Q\subset\R^d$ we write $\Gamma^\eps:=\partial G^\eps$ and $\nu_{\Gamma^\eps}$ for the outer normal vector of $G^\eps$ and consider a partial differential equation 
\begin{align}
    -\nabla\cdot\left(|\nabla u^\eps|^{p-2}\nabla \ue\right)&=f &\text{on }&Q\cap G^\eps\nonumber\\
    -|\nabla u^\eps|^{p-2}\nabla \ue\cdot\nu_{\Gamma^\eps}&=0&\text{on }&Q\cap\Gamma^\eps\label{eq:main-eps-continuous-system}\\
    \ue&=0 &\text{on }&G^\eps\cap\partial Q\,.\nonumber
\end{align}
As $\eps\to0$ we expand $\ue$ by zero to $\R^d$ and expect that $\ue\rightharpoonup u$ in $L^r(Q)$, $r\leq p$ and that $u$ solves an effective equation of the form 
\begin{align}
    -\nabla\cdot\left(\A|\nabla u|^{p-2}\nabla u\right)&=f\P(o\in G) &\text{on }&Q\\
    u&=0 &\text{on }&\partial Q\,.
\end{align}

If the above convergence behaviour holds, we say that $G$ allows homogenisation and we call $\A$ the effective conductivity of $G$. Its derivation goes far beyond simple averaging as geometric features have a major influence. As a most easy example, let $G^\eps$ be a union of finite pathwise connected components. Then the above definition of homogenisation makes no sense because $\ue$ can be shifted arbitrarily on the subsets of $G_0^\eps$ that have positive distance to $\partial Q$. From a physical point of view, with $G^\eps$ being fragmented into finite mutually disconnected sets, it is intuitive that the effective conductivity is zero.

On the other hand for many connected geometries it has been shown that, $\A>0$: Examples include periodic domains (see \cite{hopker2016extension} and references therein), minimally
smooth domains~\cite{guillen2015quasistatic} as well as the case of Bernoulli bond or site percolation~\cite{kozlov1994homogenization}.
More irregular domains have been recently investigated in \cite{heida1perforation,heida2perforation,heida3perforation} where sufficient conditions on the distributions of geometric properties such as local Lipschitz regularity or global connectivity where derived that allow to pass to the homogenisation limit.

A necessary condition for homogenisation, however, is still lacking. As a consequence, the authors proposed in \cite{stochhomIrrdomain2021} an approach for homogenisation of perforated domains where it is not clear that the perforations are ,,good enough'' but where there is also no clear indication that the domain should be too irregular for homogenisation. Such approaches using a regularisation of the random geometry can be helpful to justify a homogenised model on an irregular domain, but they leave us with a grain of salt, as it is not clear that the regularisation and the homogenisation limit really interchange.

In order to approach the question of suited or unsuited domains from the other side, in the analysis below, we will study a reasonable discrete model for a perforated domain that has the property that the effective conductivity is zero, although the microscopic geometry is topologically connected. We will also discuss heuristically which of the sufficient conditions from \cite{heida3perforation} is violated in order to make this behaviour possible, while we leave the rigorous calculation to future work. This is considered by the authors a necessary step towards more precise characterisations of admissible domains. 

Porous media and their effective conductivity $\A$ are closely related to the analysis of random walks on lattice models or in our specific case: random walks on percolation
clusters. These represent a special class of so called random conductance
models, see e.g.~\cite{MR2861133} for an extensive review. In such models, the variable-speed random walker moves along an edge at a rate equal to its conductivity. This usually admits a diffusive scaling to a Brownian motion (see e.g.~\cite{barlow2004randomwalk, andres2015degenerateEnvironment} and especially \cite{faggionato2023}) with covariance matrix $2\A$. Indeed, having zero effective conductivity
is equivalent to subdiffusivity or trapping of the related random walker. In this
regard, it is known that the random walk on the two-dimensional uniform spanning tree is subdiffusive \cite{MR2802298}
-- constituting an example of a percolating ergodic medium that features zero effective conductivity.

\begin{figure}[th]\label{fig:ERSL}
\includegraphics[width=0.5\columnwidth]{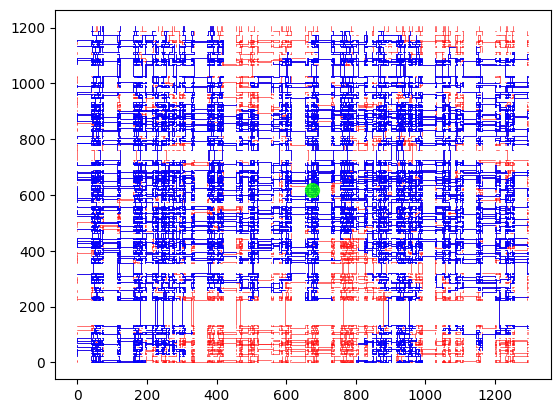}\includegraphics[width=0.5\columnwidth]{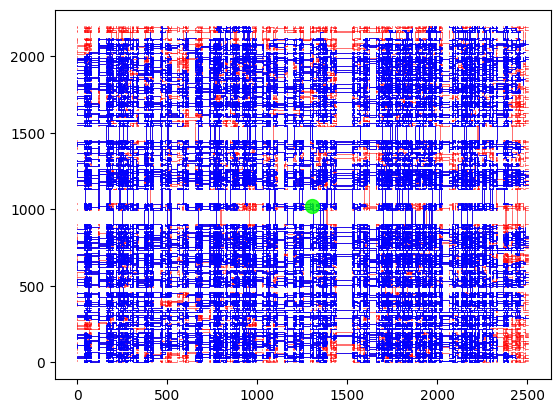}
\caption{Two realisations of the non-conductive medium, on slightly differing scales, given by the elongated randomly stretched lattice. Blue edges belong to the centrally placed
green dot's cluster (restricted to the observation window).
\label{fig:RRSL}}
\end{figure}

However, the example of a uniform spanning tree is quite artificial from a modelling perspective and very dimension dependent. Hence, in this manuscript
we present a potentially more canonical example of a non-conductive perforated medium that possesses a number of natural properties, see Figure~\ref{fig:RRSL} for illustrations. 
More precisely, the construction is based on the so-called {\em randomly stretched lattice (RSL)}, which is essentially a Bernoulli bond percolation model in a strongly correlated random environment, augmented with additional deterministic deformations.
The final model, which we call {\em elongated randomly stretched
lattice (ERSL)}, then exhibits the following key features:
\begin{enumerate}
\item It is stationary ergodic and percolates.
\item $\A=0$, i.e., the effective conductivity is zero in all directions.
\item The (annealed) probability of an edge to be open can be chosen arbitrarily
close to one.
\item The above properties can be ensured in any dimension larger or equal to two.
\end{enumerate}

Let us comment on some potentially simpler models that however violate at least one of the above conditions. 
First, the planar uniform spanning tree satisfies the first two conditions, but
it is fractal and it is not connected in high dimensions. Second, it is not hard to
construct a percolating lattice model that is non-conductive in all directions
except one, see e.g.~Section~\ref{sec:outlook} for a brief description.

Let us finally mention that, while the ERSL is presented and while all relevant calculations are done in the discrete
setting, they can be easily brought into the continuum: one simply needs to thicken the ERSL by some $r<1/2$ after embedding the lattice into $\R^{d}$.

\section{Setting and main result}

The elongated randomly stretched lattice (ERSL) is (as the name suggests) a translation- and (lattice)-rotation-invariant nearest-neighbour percolation model on $\Z^d$, constructed by elongating a randomly stretched lattice (RSL). We define both models in Section~\ref{sec_RSL}. Given unit conductance on open edges, let $\A$ be the effective conductivity of the ERSL. As mentioned in the introduction, while $\A$ is a semi-positive definite symmetric matrix, $2\A$ is also the covariance matrix of the diffusively scaled (variable speed) random walker on the ERSL. Moreover, all symmetries of the ERSL transfer to $\A$, and in particular $\A$ obeys lattice-rotation invariance. 

In order to specify the connection between $\A$ and the ERSL, let us write $[0,n]:=\{0,1,\dots,n\}$ and $y\sim x$ if and only if there exists an
open edge between vertices $x$ and $y$ in the ERSL. The effective conductivity in the direction $\e_1$ is represented by the following asymptotic minimisation problem, see Section \ref{sub:effective-conductivity}:
\begin{equation}
\e_{1}^{\mathrm{t}}\A\e_{1}=\lim_{n\to\infty}n^{2-d}\inf_{V\in\mathcal{D}_{n}}\tfrac{1}{2}\sum_{\substack{z,\Tilde{z}\in[0,n]^{d}\\
z\sim \Tilde{z}
}
}|V(\Tilde{z})-V(z)|^{2},
\label{eq:eff-cond-discrete}
\end{equation}
which holds for almost-every realisation of the ERSL, due to ergodicity. Here, $\e_1$ denotes the unit vector in the first coordinate with $\e^{\mathrm t}_1$ its transposition and $\mathcal{D}_{n}$
consists of functions $V:[0,n]^{d}\to\R$ satisfying
\[
V(0,z_{2},\dots,z_{d})=0\quad\text{and}\quad V(n,z_{2},\dots,z_{d})=1\,.
\]
Under the lattice isotropy, we have $\e_i^{\mathrm{t}}\A\e_i=\e_{1}^{\mathrm{t}}\A\e_{1}$ for all $i\le d$. In fact, since $\A$ is symmetric and rotationally invariant under lattice rotations, we even have isotropy of $\A$, i.e.~$\A=a_0 I_d$ for some $a_0\in\R$ and identity matrix $I_d$.
Let us mention that, by the Dirichlet principle, any minimiser $V_{\mathrm{min}}$ is harmonic, i.e.~it
even satisfies
\[
\sum_{\Tilde{z}\sim z}(V_{\mathrm{min}}(\Tilde{z})-V_{\mathrm{min}}(z))=0\quad\text{ for all }z\in(0,n)\times[0,n]^{d-1},
\]
see e.g.~\cite{DoyleSnell}. 
Before we define the ERSL precisely let us state our main result.
\begin{thm}\label{thm:main-thm}
For any $d\geq2$ and $\bar{p}\in(0,1)$, there exists an ERSL -- a stationary
ergodic nearest-neighbour bond percolation model on $\Z^{d}$ -- satisfying
the following properties.
\begin{enumerate}
\item The ERSL percolates almost surely and is lattice rotation invariant.
\item $\P(e\text{ is open in the ERSL})\geq\bar{p}$ for any edge $e$ in
$\Z^{d}$.
\item For the associated conductivity $\A$, as defined in~\eqref{eq:eff-cond-discrete}, we have that $\A=0$.
\end{enumerate}
In particular, the random walk on the ERSL is subdiffusive.
\end{thm}

In the following construction as well as all proofs, we will restrict ourselves to the planar case, $d=2$, for convenience. All other cases $d>2$ follow by completely analogous arguments.

\section{Construction: RSL and ERSL}\label{sec_RSL}

Let us first introduce a prototypical lattice model with columnar
disorder: the randomly stretched lattice. It is a bond percolation
model on $\Z^{2}$ where entire columns are made ``weak'', i.e., bonds
in such areas are likely to be closed.
\begin{defn}[Randomly stretched lattice (RSL)]
\label{def:RSL}
Let $p,q\in(0,1)$ and consider families $N^{(x)}:=(N_{i}^{(x)})_{i\in\Z}$
and $N^{(y)}:=(N_{j}^{(y)})_{j\in\Z}$ of iid geometric random variables satisfying
\[
\P(N_{0}^{(x)}\geq \ell+1)=\P(N_{0}^{(y)}\geq \ell+1):=q^{\ell}\,.
\]
 Given a realisation of $N^{(x)}$ and $N^{(y)}$, all the bonds in
$\Z^{2}$ are open independently with probabilities
\[
\P\Big((i,j)\leftrightarrow(i+1,j)\text{ is open}\,\vert\,N^{(x)},\,N^{(y)}\Big):=p^{N_{i}^{(x)}}
\]
and
\[
\P\Big((i,j)\leftrightarrow(i,j+1)\text{ is open}\,\vert\,N^{(x)},\,N^{(y)}\Big):=p^{N_{j}^{(y)}}\,.
\]
This model is called the {\em randomly stretched lattice (RSL)}.
We will often refer to the value $\ell$ as the {\em badness}. 
\end{defn}
The RSL features a non-trivial percolation behaviour in the sense that, in non-trivial parameter regimes, realisation of the RSL contain unbounded connected components with positive probability, and in other non-trivial regimes not. More precisely, we have the following statement. 
\begin{thm}[Existence of supercritical regime in the RSL, \cite{MR1761579,MR2116736,MR4634238,rsl2023}]
\label{thm:RSL-Supercritical}~Consider the RSL as in Definition~\ref{def:RSL}
with $p>1/2$. Then, there exists $q_{c}\in(0,1)$ such that,
for all $q\leq q_{c}$, the RSL percolates almost-surely. 
\end{thm}

\begin{proof}
For $d\geq3$, percolation has been shown in \cite{MR1761579}, while
the $d=2$ case was established in \cite{MR2116736} for large $p$.
This result as well as methods have been improved over time in \cite{MR4634238,rsl2023}.
\end{proof}
Let us mention that, in two dimensions, we may even ensure finite (albeit not uniformly bounded) size void spaces using circuits of open bonds
around $\Lambda$ for every finite set $\Lambda\subset\Z^{2}$. This is shown in~\cite{MR4660696} for a part of the supercritical regime, however, the approach in~\cite{rsl2023} enables a much simpler proof covering all $p>1/2$ (not necessarily the whole supercritical regime though).

\medskip
Let us highlight that, while the RSL features infinitely long dependencies, these dependencies are confined to columns and rows. Therefore, the RSL is mixing in all diagonal directions which yields ergodicity.

\medskip{}

\noindent\fbox{\begin{minipage}[t]{1\columnwidth - 2\fboxsep - 2\fboxrule}%
In the following, we will fix a parameter pair $p,q\in(0,1)$ for
which the RSL percolates and additionally $p>q$. Furthermore, whenever we refer to the RSL in the future, we mean a
realisation. All statements relating to the RSL are meant in the almost-sure sense.
\end{minipage}}

\medskip{}

Unfortunately, we are unable to establish $\A=0$ directly for the RSL. While weak columns have high resistance,
they do not occur frequently enough. Fortunately, duplicating (or rather deterministically elongating) columns and rows solves the issue. Doing so has no impact on the connectivity of the underlying percolation model, but it has a huge effect on the conductivity. Roughly speaking, the elongation is done such that ``bad layers'' become exponentially large.  

We describe this procedure in the following only for $(N_{i})_{i\in\Z}:=(N_{i}^{(x)})_{i\in\Z}$, i.e., along the first coordinate. All other coordinates are treated in the same fashion. Recall that $\lceil a\rceil:=\inf\{n\in \Z\colon a\le n\}$. The key idea is to use the labels $N_i$ in the RSL to insert additional columns into the lattice. More precisely, let $\sigma\in(0,1)$ such that $q^{\sigma}>p$ and set
\begin{align}\label{eq_elong}
\S(\ell):=\lceil q^{-\ell(1-\sigma)}\rceil,\qquad l\in\N.  
\end{align}
Given the $i$-th column in the RSL, with label $N_{i}$, we deterministically elongate this column to have width $\S(N_i)$. In other words, for any realisation of the RSL, for any column $i$, we will insert additionally $\S(N_i)-1$ copies of said column's  horizontal edges including their open or closed state. We will call this elongated strip of
width $\S(\ell)$ a {\em layer of badness} $\ell$. It consists of $\S(\ell)$ many {\em columns of badness $\ell$}. 

\begin{figure}[th]
\includegraphics[width=0.33\columnwidth]{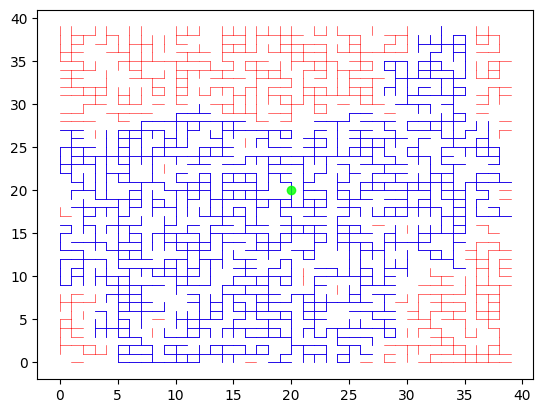}\includegraphics[width=0.33\columnwidth]{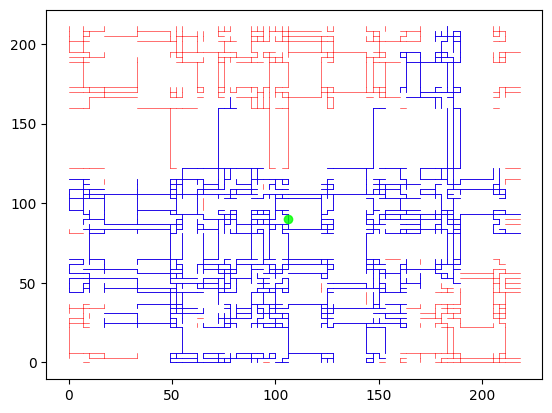}
\includegraphics[width=0.33\columnwidth]{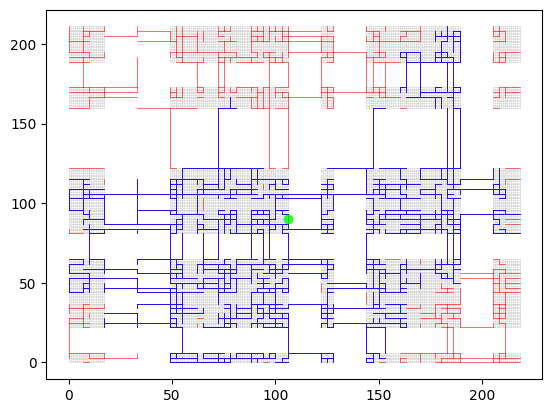}\caption{Realisations of the RSL with parameters $p=0.65,q=0.3$ (left), the elongated version with $\sigma=0.25$ (middle), and the filled version $F_L(\text{RSL})$ with $L=2$ (right). 
Blue edges belong to the centrally placed green dot's connected component (before using the grey filling).\label{fig:elongation}}
\end{figure}

Concerning the openness or closedness of the edges in the layers we employ an additional modification using the parameter $L\in \N$, to be specified later. If an edges lies in a rectangle spanned by a horizontal and a vertical layer with badness $l\leq L$ then, then we set it to be open. Otherwise, they remain unmodified. Recall that we do the same analogously for rows, leading to the {\em elongation and filling transformation} $\text{RSL}\mapsto F_L(\text{RSL})$, as illustrated in Figure~\ref{fig:elongation}.

\medskip
We have to be cautious as the deterministic elongation destroys the stationarity of the lattice. But we can introduce a stationarising random initial shift
$\bar Z=(Z^{(x)},Z^{(y)})$ (with $Z^{(x)},Z^{(y)}>0$ iid for both coordinates) to recover stationarity, see e.g.~\cite[Theorem 9.1]{Thorisson}. This can be done as long as the expected elongation is finite, which is the case since, by the definition in~\eqref{eq_elong}, 
$$\E[\S(N_1)]=(1-q)\sum_{\ell\ge 1}q^\ell\lceil q^{-\ell\cdot(1-\sigma)}\rceil\leq 1 + \frac{1-q}{1-q^\sigma}<\infty.$$
To be precise, $Z^{(x)}$ has the size-biased distribution $$\P(Z^{(x)}= n):=n\E[\S(N_1)]^{-1}\P(\S(N_1)= n)$$ 
and we note that $Z^{(x)}$ is finite but possibly without finite first moment. The probability of a (horizontal) edge in the initial layer of badness $Z^{(x)}$ to be open is given by $p^{N'_0}$, where $N'_0$ is a random variable with $$\P(N'_0= n)=\S(n)\E[\S(N_1)]^{-1}\P(N_1=n).$$
Now we can define our model. 
\begin{defn}[Elongated randomly stretched lattice (ERSL)]
\label{def:ERSL}
The elongated randomly stretched lattice (ERSL) is defined as $F_L(\text{RSL}')-\bar U$, where $\text{RSL}'$ is the RSL in which the initial marks $N^{(x)}_0$ and $N^{(y)}_0$ are replaced by iid copies of the size-biased marks $N'_0$ as defined above. The random shift $\bar U=(U^{(x)},U^{(y)})$ is given by independent uniformly distributed random variables $U^{(x)}\in [0,\S(N^{(x)}_0)]$.  
\end{defn}
Let us mention that it is precisely this possibility to create a shift-invariant version that fails for an elongated, non-conductive version of simple Bernoulli bond percolation on $\Z^2$. Indeed, to get infinite resistance, the elongations would need to have infinite first moment. 

Let us collect the first properties of the ERSL.
\begin{proof}[Proof of Theorem~\ref{thm:main-thm} Part (1) and (2)]
Note that the ERSL is ergodic since it is still mixing along the diagonal directions. Furthermore, it is a nearest-neighbour bond percolation model and its distribution is invariant under all lattice rotations. In the given parameter regime $p>1/2$ and $q\leq \min\{q_c,p\}$ it features percolation almost surely and this percolation is maintained if we set additional edges to be open under the $L$-dependent rule. In particular, for all $\bar p$, we can pick $L$ sufficiently large such that the typical edge in the system is open with probability not smaller than $\bar p$. 
\end{proof}

The next section verifies that $\A=0$.

\section{Checking non-conductance }
Before we start, let us give some intuition. Disregarding the stationarisation as well as the $L$-dependent rule, the RSL with both dimensions being elongated serves as our material model in $\Z^{2}$. The idea behind the introduction of the parameter $\sigma$ is now the following: In macroscopic columns of badness $\ell$, the fraction of open edges is at most $p^\ell$. This means, that the whole column behaves similarly to one where all edges are open but with conductance $p^\ell$ instead of unit conductance. In other words, the resistance is not smaller than approximately $p^{-\ell}$. Additionally, the probability of seeing a column of badness $\ell$ is given by $\S(\ell)q^\ell (1-q) \E[\S(N_1)]^{-1} \approx Cq^{\sigma \ell}$
instead of $q^{\ell}(1-q)$ due to size bias, for some constant $C>0$. Observe now that $\sigma$ has been chosen such that $q^\sigma > p$ and hence, the expected macroscopic resistance is given by
$$C\sum_{\ell=1}^\infty  p^{-\ell}q^{\sigma \ell}=C\sum_{\ell=1}^\infty  (q^{\sigma}/p)^\ell = \infty.$$
But this means that we have zero conductivity. 

On a technical level, it suffices to show that the right-hand side in~\eqref{eq:eff-cond-discrete} equals $0$. 
For this, the idea is to locate bad layers and ensure that only few edges are open. In order to do that, let us introduce some relevant quantities that will aid
us in our calculations.

\subsection{Parameters}
Given $p,q,\sigma$ as before, we also consider
$$ \gamma:=\Big[(1-\sigma)+\frac{\log(p)}{\log(q)}\Big]^{-1/2}.$$
As $q^\sigma > p$, i.e., $\sigma<\frac{\log(p)}{\log(q)}$, we see that $\gamma\in (0,1)$. We will see that we find layers with badness $\ell_{n}$ in boxes of size $n$ with high probability, where
\[
\ell_{n}:=\frac{\gamma\log n}{\log(q^{-1})}\,.
\]
Let us note that, with these quantities, we have $p^{\ell_{n}}\geq q^{\ell_{n}}n\to\infty$ and, as $\gamma<1$, also
\begin{align*}
   \frac{p^{\ell_{n}}n}{\S(\ell_{n})} &\leq \exp \Big(\log(n)\Big[\frac{\gamma}{\log(q^{-1})}\log(p) + 1 + \frac{\gamma}{\log(q^{-1})}\log(q)(1-\sigma)\Big]\Big) \\
&= \exp \Big(\log(n)\Big[-\gamma\Big(\frac{\log(p)}{\log(q)}+(1-\sigma)\Big) + 1\Big]\Big) = n^{1-1/\gamma}\to 0. 
\end{align*}

\subsection{Finding bad layers}

We want to estimate the probability of finding suitably weak columns as $n\to\infty$. Due to the stationarising (horizontal) shift $U=U^{(x)}$, we know that $[U,U+\S(N_1))$ is a layer in $[0,\infty)$. Analogously, the $i$-th consequent layer is $[U+\sum_{k=1}^{i}\S(N_{k}),U+\sum_{k=0}^{i+1}\S(N_{k}))$
and has badness $N_{i+1}$. Therefore, finding a layer with badness not smaller than $\ell_{n}$ inside $[0,n]$ is guaranteed under the event
\[
E_{n}:=\Big\{\exists i\in\N:\ N_{i}\geq \ell_{n}\quad\text{and}\quad U+\sum_{k=1}^{i}\S(N_{k})<n\Big\}.
\]

\begin{lem}[Probability of bad layers]\label{lem_1}
We have that $
\lim_{n\uparrow\infty}\P(E_{n})=1$.
\end{lem}

\begin{proof}
Set $C:=(2\E[\S(N_{1})])^{-1}$ and denote the events
\begin{align*}
A_{n} & :=\Big\{\sum_{k=1}^{Cn}\S(N_{k})<n-U\Big\}=\Big\{\frac{1}{Cn}\sum_{k=1}^{Cn}\S(N_{k})<2\E[\S(N_{1})]-\frac{1}{Cn}U\Big\}\text{ and}\\
B_{n} & :=\big\{\exists 1\leq i\leq Cn:\ N_{i}\geq \ell_{n}\big\}.
\end{align*}
 Then, by restricting to $i\leq Cn$, we see that $E_{n}\supset A_{n}\cap B_{n}$.
But, by the law of large numbers and almost-sure finiteness of $U$, we have $\P(A_{n})\uparrow 1$. Further, since the $N_{i}$ are iid geometric random variables, we have
\[
1-\P(B_{n})=\big(1-q^{\ell_{n}}\big)^{Cn}\leq \exp(-Cq^{\ell_{n}}n)\xrightarrow{n\to\infty}0\,,
\]
by the choice of $\ell_{n}$. Combining the two statements yields
\[
\P(E_{n})\ge \P(A_{n}\cap B_{n})\ge \P(A_{n})-\P(B^c_{n})\xrightarrow{n\to\infty}1.
\]
which finishes the proof.
\end{proof}
Each column inside the box $[0,n]^{2}$ contains at most $n$ horizontal edges. Now, we check that only few of these edges are open with high probability if they belong to a bad layer.

\begin{lem}[Conductivity through bad layers]\label{lem_2}
Assume that $E_{n}$ occurs and consider one associated layer
of badness at least $\ell_{n}$. Let $F_{n}$ be the event that at
most $2np^{\ell_{n}}$ out of at most $n$ horizontal edges inside
the layer are open. Then,
\[
\lim_{n\uparrow \infty}\P(F_{n}\cap E_n)=1.
\]
\end{lem}

\begin{proof}
Let $X_{1},\dots,X_{n}$ be iid Bernoulli random variables with $\P(X_{1}=1)=p^{\ell_{n}}$
and $Y:=\sum_{i=1}^{n}X_{i}$. Then, by the Chebyshev inequality,
\begin{align*}
\P(F_n^c\cap E_n)&\le \P(Y-\E[Y]>\E[Y])\leq\frac{n\mathrm{Var}[X_{1}]}{(n\E[X_{1}])^{2}}=\frac{1}{n\E[X_{1}]}\cdot\frac{\mathrm{Var}[X_{1}]}{\E[X_{1}]}\leq\frac{1}{p^{\ell_{n}}n}\cdot1 \xrightarrow{n\to\infty} 0,
\end{align*}
which shows the claim by Lemma~\ref{lem_1} and the choice of $\ell_{n}$.
\end{proof}

\subsection{Calculating the effective conductivity}

We may finally calculate the right-hand side in~\eqref{eq:eff-cond-discrete}. 
\begin{lem}[Upper bound]\label{lem_3}
Under the event $E_n\cap F_n$ there exists a $V\in\mathcal{D}_{n}$ such that 
\begin{equation}\label{eq_2}
\sum_{z,\Tilde{z}\in[0,n]^{2}\colon 
z\sim \Tilde{z}
}|V(\Tilde{z})-V(z)|^{2}\leq4n^{1-1/\gamma}.
\end{equation}
\end{lem}

\begin{proof}
Under the event $F_n\cap E_n$ there exists a layer with badness at least $\ell_{n}$
inside $[0,n]$. Let $X$ denote its starting location. Then, we define $V:[0,n]^{2}\to\R$ as 
\[
V(i,j):=\begin{cases}
0 & i<X,\\
1 & i>X+\S(\ell_{n}),\\
\frac{i-X}{\S(\ell_{n})} & i\in[X,X+\S(\ell_{n})],
\end{cases}
\]
and note that $V(i,j)$ does not depend on $j$ and is constant except inside the bad layer where it linearly grows to $1$. Therefore, the only contribution to the sum in~\eqref{eq_2}
comes from horizontal edges involving $i\in[X,X+\S(\ell_{n})]$ and in particular $V\in\mathcal{D}_{n}$.

Let us focus on the chosen bad layer. There, under $F_n$, we have at most $2np^{\ell_{n}}$ open edges along a strip of size $\S(\ell_{n})$ (which is part of a potentially larger layer), so the contribution is
\begin{align*}
\sum_{\substack{z,\Tilde{z}\in[0,n]^{2}\\
z\sim \Tilde{z}
}
}|V(\Tilde{z})-V(z)|^{2} & \leq 4np^{\ell_{n}}\S(\ell_{n})S(\ell_{n})^{-2}\leq  4n^{1-1/\gamma},
\end{align*}
as desired.
\end{proof}

\begin{proof}[Proof of Theorem~\ref{thm:main-thm} Part (3)]
Note that by the Borel--Cantelli lemma, there exists a subsequence $(n_k)_{k\ge 1}$ such that almost surely, for all but finitely many $k$, the event $E_{n_k}\cap F_{n_k}$ occurs. Hence, using Lemma~\ref{lem_3} yields
\begin{align*}
\P\Big(\lim_{n\uparrow\infty}&\inf_{W\in\mathcal{D}_{n}}\sum_{x,y\in[0,n]^{d}\colon 
x\sim y
}|W(y)-W(x)|^{2}\le 0\Big)\\
&= 
\P\Big(\inf_{W\in\mathcal{D}_{n_k}}\sum_{x,y\in[0,n_k]^{d}\colon 
x\sim y
}|W(y)-W(x)|^{2}\le 4n_k^{1-1/\gamma}\text{ for infinitely many } k\Big)\\
&\geq \P\Big(E_{n_k}\cap F_{n_k}\text{ happens for infinitely many }k\Big)=1\,,
\end{align*}
where the first equality follows from the fact that the limit in Equation \eqref{eq:eff-cond-discrete} exists almost surely as well as $\lim_{k\to\infty}n_k^{1-1/\gamma} = 0$.
This shows $\A=0$.
\end{proof}

\section{Background on discrete models for perforated domains}

\subsection{Justification of discrete models replacing continuous problems}\label{sub:justification-discrete}
We will now demonstrate that discrete homogenisation problems can at least in some cases be considered as continuous homogenisation problems. The implication of this insight is that, constructing a discrete medium that is topologically connected but has macroscopic conductivity zero can be mapped onto a continuous medium with the same properties.

The classical point of view, which we discuss first, somehow follows the opposite direction. However, we provide it here for completeness.

\paragraph{The classical point of view.}
Historically, the upscaling of discrete models were first proposed independently in \cite{kunnemann1983diffusion,kozlov1987averaging} as substitutes for the homogenisation of partial differential equations.
A basic idea behind this discretisation is the finite-volume approach: A partial differential equation of the form $-\nabla\cdot(a(x)\nabla u)=f(x)$ can be discretised on a cubic grid by $-\delta^{-2}\sum_\pm\sum_{j=1}^d a_{x,j}\left( u(x\pm\delta e_j)-u(x)\right)=f(x)$, where $x\in\delta\Z^d$ and $a_{x,j}$ is constructed properly in~\cite{droniou2018gradient}.

Hence, if we consider a homogenisation problem $-\nabla\cdot(a(x/\eps)\nabla u)=f(x)$, we can consider instead 
$-(\delta\eps)^{-2}\sum_{j=1}^d a_{x,j} \left(u(x+\delta\eps e_j)-u(x)\right)=f(x)$, where $x\in\eps\delta\Z^d$ or equivalently, after a rescaling,
$$-\eps^{-2}\sum_{j=1}^d a^\eps_{i,j} \left(u(x_i+\eps e_j)-u(x_i)\right)=f(x_i)\,,\qquad x\in\eps\Z^d\,.$$

In order to transfer this insight to the case of a perforated domain, we can consider a stationary ergodic random domain with holes that are large compared to the grid distance in $\Z^d$. Then, we consider $G_\Z:=G\cap\Z^d$ and say $a_{i,j}=1$ if and only if $x_i,\,x_i+e_j\in G_\Z$ and $a_{i,j}=0$ otherwise. This mimics the behaviour of \eqref{eq:main-eps-continuous-system} in the discrete setting.

\paragraph{The exact solution point of view.}
Through another point of view, our discrete solutions can be mapped one on one to solutions for a subclass of problems on a special perforated domain. 
In order to avoid struggles with boundary conditions, we consider the full-space problem, even though the major events happen around a bounded domain $Q$.

For every vertex $x$ of our rectangular grid, we consider for $\delta\ll 1$ the cube of width $\delta$ with centre $x$ and call it $\C(x,\delta)$. If $x\sim y$ are connected neighbours in our discrete model, we connect the two cubes $\C(x,\delta)$ and $\C(y,\delta)$ by their combined convex hull $\C(x,y,\delta)$, e.g., a rectangular cylinder with a $(d-1)$-dimensional cube of size $\delta$ as its base. We call $\C^0(x,y,\delta)=\C(x,y,\delta)\setminus(\C(x,\delta)\cup\C(y,\delta))$ and $G^\delta=\bigcup_x \C(x,\delta)\cup\bigcup_{x\sim y}\C^0(x,y,\delta)$. 

Given values $f_x$ for each vertex $x$, with $f_x=0$ for $x\not\in Q$ we consider the discrete equation for $u$
\begin{equation}\label{eq:discr-graph-version-1}
\forall x\in\Z^d\cap Q\,,\quad \sum_{x,y\in[0,n]^{d}\colon
x\sim y}  u_x-u_y=f_x\,\qquad\text{ and}\qquad \forall x\Z^d\setminus Q,\quad u_x=0. 
\end{equation}
This problem has a unique solution as the linear map on the left-hand side is positive definite.

Next we define $f^\delta(x)=f_x$ on $\C(x,\delta)$ and $f^\delta(x)=0$ else. 
We set $u^\delta$ as the linear interpolation of $u_x$ on $\partial\C(x,\delta)\cap\partial\C^0(x,y,\delta)$ and $u_y$ on $\partial\C(y,\delta)\cap\partial\C^0 (x,y,\delta)$. 
Furthermore, let $u^\delta$ solve $-\Delta u^\delta=f$ on $\C(x,\delta)$ with $u^\delta=u_x$ on $\partial\C(x,\delta)$. 
Then, $u^\delta$ is an $H^1_{loc}(G^\delta)$ function that solves $-\Delta u^\delta=f^\delta$ on $G^\delta$. 
Furthermore, $u^\delta=0$ outside a sufficiently large ball around $Q$.

It is thus reasonable to consider a sequence of discrete solutions as a sequence of solutions to~\eqref{eq:main-eps-continuous-system}, turning the discrete homogenisation problem into a continuous homogenisation problem.

\subsection{Formulas for the effective conductivity}\label{sub:effective-conductivity}

We will now justify our formula for the effective conductivity~\eqref{eq:eff-cond-discrete}. Since this formula is well established in literature, it is not our goal to rigorously derive it, but to recap some of the main arguments as to why this formula is correct.
Given a lattice $\L\subset\Z^d$, we consider $\L^\eps :=\eps\L$ as well as the following scaled version of \eqref{eq:discr-graph-version-1}
\begin{equation}\label{eq:discr-graph-version-2}
\forall z\in\L^\eps\cap Q\,,\qquad \eps^{-2}\sum_{\substack{z,\Tilde{z}\in[0,n]^{d}\\
z\sim \Tilde{z}}}  u_z^\eps-u_{\Tilde{z}}^\eps=f_z\, ,    
\end{equation}
which takes the following form by a variational principle
$$\ue = \arg\min u\mapsto \eps^{d}\sum_{z\in\L^\eps\cap Q}\left(\frac12\eps^{-2}\sum_{\Tilde{z} \sim z}|u_z-u_{\Tilde{z}}|^2-f_z u_z\right).$$
Boundary conditions can be imposed by restricting the space over which the minimum is taken.

Without going into detail, but referring to~\cite{MR2861133}, the effective conductivity is again defined as $\A$ such that $f^\eps\to f$ and $\ue\to u$ in an appropriate sense (this involves mapping discrete functions to continuous ones) implies that $u$ solves $-\nabla\cdot(\A\nabla u)=f$ on $Q$. 
Furthermore, using $\Gamma$-convergence arguments, one can draw the conclusion that the minimisers $\ue$ from above satisfy
\begin{equation}\label{eq:conv-mini-gamma}
\int_Q \nabla u\cdot\A\nabla u- fu{\rm d} x = \liminf_{\eps\to0}\eps^{d}\sum_{z\in\L^\eps\cap Q}\left(\frac12\eps^{-2}\sum_{\Tilde{z}\sim z}|u_z-u_{\Tilde{z}} |^2-f_z u_z\right)\,.    
\end{equation}

With regard to \eqref{eq:eff-cond-discrete}, let us note that in the continuous case, a function satisfying $u(0,z_2,\dots,z_d)=0$, $u(1,z_2,\dots,z_d)=1$,  $-\nabla\cdot(\A\nabla u)=0$, and also minimising the left-hand side of \eqref{eq:conv-mini-gamma}, has to satisfy $u(z)=z_1$ and it holds that $\int_{(0,1)^d}\nabla u\cdot\A\nabla u{\rm d} x=\e_1\A \e_1$. 

The correctness of \eqref{eq:eff-cond-discrete} now follows from a rescaling, choosing $\eps=N^{-1}$.

\subsection{Violation of the homogenisation conditions in \cite{heida3perforation}}

In~\cite{heida3perforation} it is assumed that we can distribute a point process $X=(X_i)_{i\ge 1}$, inside the random geometry, that is jointly stationary and with a uniform minimal distance $\delta>0$ to the boundary. Then, these points $X_i$ are used to create a Voronoi tessellation where each cell $C_i$, corresponding to $X_i$, has a diameter $D_i$. 
Within our above construction of a channel network, this situation can be reproduced for example by choosing a subset of $\Z^d\cap G^\delta$, as each of these points has a distance $\delta$ to $\partial G^\delta$. 

Now,~\cite{heida3perforation} states three conditions on the random geometry and the chosen point process that together ensure positive conductivity. 
Two of these conditions are concerned with the moments of local Lipschitz regularity and thickness of pipes, which are both satisfied even uniformly in our model. 
The third condition, \cite[Equation (1.12)]{heida3perforation}, is related to the Voronoi cells emerging from $X$ and it implies that the $(5d+1)$-th moment of the typical diameter exists, i.e., $\E[D_0^{5d+1}]<\infty$, where $D_0$ is the diameter corresponding to a uniformly chosen point in the domain, see Palm theory for details.   
Since we characterise our geometry by exponential distributions and the conditions in~\cite{heida3perforation} are polynomial, let us shortly sketch how we believe this fits together, i.e., why our model violates the sufficient conditions in~\cite{heida3perforation}.

We observe that the diameters $D_i$ are related to the thickness of bad layers and to the mean mutual distance of channels in bad layers. 
Concerning the first part, we note that $\P(\S(N_0)=q^{-l(1-\sigma)})\approx q^l$ and thus $\P(k<D_0\le k+1)\sim k^{1/(\sigma-1)}$. Concerning the second part, we can put ourselves in one of the channels in the layer of badness $l$ and observe that the probability to find another channel in orthogonal direction to the current channel and within a distance $R$ is proportional to $\omega_l (1-p^l)^{R^{d-1}}$ where $\omega_l$ is supposedly exponentially decreasing. Comparing the sum $\sum_l \omega_l(1-p^l)^{R^{d-1}}$ with an integral $\int_1^\infty \omega_l(1-p^x)^{R^{d-1}}\mathrm{d}x$, the probability becomes polynomial in $R^{1-d}$. 

We leave the detailed verification of the above heuristic to future investigations as we expect the corresponding calculations to be involved, for example due to boundary effects in the creation of the Voronoi cells. 

\section{Discussions and outlook}\label{sec:outlook}
\begin{figure}[th]
\includegraphics[width=1\columnwidth]{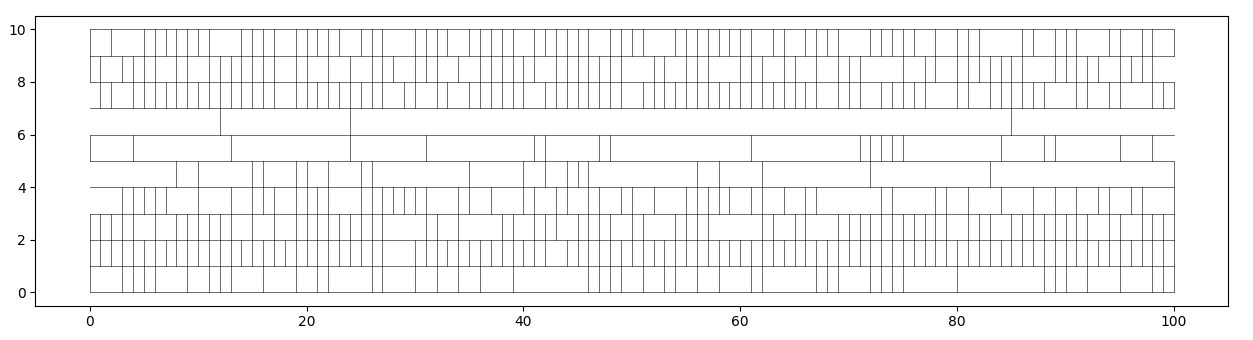}
\caption{Realisation of a lattice model with zero conductivity in vertical direction. Here, $P_0$ is a uniform random variable.}\label{fig:one-directional}
\end{figure}

We briefly mentioned a simple percolation model with infinite resistance in all (standard lattice) directions but one, which we illustrate in Figure~\ref{fig:one-directional}: Let us first choose $\e_1\in\Z^d$ as our special direction. In this direction, we set all edges $(v\leftrightarrow v+\e_1)$ to be open. For all other directions, we do the following (only illustrated for $d=2$ and analogously for higher dimensions): For each column $i$ of horizontal edges, sample independently a $P_i\in(0,1)$. Then, the edges in said column are independently set open with probability $P_i$ and closed otherwise. If all the $P_i$ are iid and $\E[P_0^{-1}]=\infty$, then the expected resistance is infinite and the effective conductivity is zero in this direction. In fact, this model features exactly one connected component which contains all open edges.

Previously, we exploited the existence of weak columns in the RSL by elongating them. With actual control on the percolation cluster, one might be able to answer the following question: Is there perhaps a supercritical regime in the unmodified RSL which already features zero conductivity on large scales?

By now, we have shown that the (variable speed) random walker on a properly tuned ERSL is subdiffusive. Figuring out the exact scaling and determining whether it depends on the elongation parameter $\sigma$ seems worthwhile.

\begin{acknowledgement*}
We would like to thank Sebastian Andres, Marek Biskup and Alessandra Faggionato for inspiring discussions. BJ and ADV received support by the Leibniz Association within the
Leibniz Junior Research Group on \emph{Probabilistic Methods for Dynamic Communication Networks} as part of the Leibniz Competition (grant no.\textbackslash{} J105/2020). MH was funded by DFG through SPP2256 project HE 8716/1-1 project ID: 441154659.
\end{acknowledgement*}

\printbibliography[heading=bibintoc]

\end{document}